\documentclass[11pt]{article}
\usepackage[margin=25mm]{geometry}
\usepackage{amsmath,amssymb}
\usepackage{microtype}
\usepackage{hyperref}

\title{A Hadamard Cube of Order 22}
\author{Amin Bahmanian\\Department of Mathematics, Illinois State University\and
Sho Suda\\Department of Mathematics, National Defense Academy of Japan}
\date{October 1, 2026}

\begin{document}
\maketitle

\begin{abstract}
We give an explicit three-dimensional Hadamard matrix of order $22$.
The example is presented by listing its twenty-two $22\times22$ two-dimensional slices.
\end{abstract}

\section{Introduction}
Let $[v]=\{1,\dots,v\}$. A three-dimensional Hadamard matrix, or Hadamard cube, of order $v$ is a map
\[
 H:[v]^3\longrightarrow\{\pm1\}
\]
such that, for every $j\in\{1,2,3\}$ and every $a,b\in[v]$,
\[
 \sum_{1\le i_1,\ldots,\widehat{i_j},\ldots,i_3\le v}
 H(i_1,\ldots,a,\ldots,i_3)
 H(i_1,\ldots,b,\ldots,i_3)
 =v^2\delta_{ab}.
\]
Equivalently, the $v$ two-dimensional layers perpendicular to each of the three coordinate directions are mutually orthogonal.

For $\ell\in\{1,2\}$, an $\ell$-layer of a three-dimensional $\{\pm1\}$-array is a subarray obtained by fixing $3-\ell$ coordinates. An $(v,3,\ell)$-Hadamard hypercube is a three-dimensional $\{\pm1\}$-array of order $v$ in which every two parallel $\ell$-layers are orthogonal. The integer $\ell$ is called the \emph{class} of the Hadamard hypercube. Thus, a Hadamard cube as defined above is an $(v,3,2)$-Hadamard hypercube.

Kr\v{c}adinac, Pav\v{c}evi\'c, and Tabak constructed three-dimensional Hadamard matrices of many previously unknown orders and identified $22$ as the smallest order congruent to $2$ modulo $4$ not covered by their construction; they explicitly stated that the existence of a three-dimensional Hadamard matrix of order $22$ was open~\cite{KPT2023}. The case also appeared among the open problems collected in~\cite{KZXS2025}.

More recently, Bahmanian and Suda introduced two further constructions of Hadamard hypercubes, one based on conference matrices and one recursive, and thereby resolved many previously open cases~\cite{BS2026}. These constructions still do not provide an example of order $22$. The purpose of this note is simply to record an explicit Hadamard cube of order $22$.

\section{An explicit example}
Let $H=(H(i,j,k))_{1\le i,j,k\le22}$. For $i=1,\dots,22$, define the slice
\[
 H_i=(H(i,j,k))_{1\le j,k\le22}.
\]
The twenty-two slices are displayed below. We write $+$ and $-$ for $+1$ and $-1$, respectively.

\clearpage
\newgeometry{left=12mm,right=12mm,top=10mm,bottom=12mm}
\begin{minipage}[t]{0.485\textwidth}
\centering
{\footnotesize $H_{1}=$}\par\vspace{-0.2mm}
{\fontsize{5.8}{5.9}\selectfont
\setlength{\arraycolsep}{0.55pt}\renewcommand{\arraystretch}{0.68}
\(\left(\begin{array}{*{22}{c}}
- & - & - & - & + & - & - & - & - & - & - & - & + & - & + & - & + & - & - & - & + & - \\
- & + & - & + & + & - & + & - & - & - & + & + & - & + & - & - & + & - & - & + & + & + \\
- & - & - & - & + & + & + & + & + & - & + & + & - & + & - & + & + & + & + & - & + & - \\
- & - & - & + & + & - & - & - & + & + & - & + & + & + & - & - & + & + & + & + & - & - \\
+ & + & + & + & + & - & + & - & + & - & + & + & + & - & - & + & - & - & - & + & - & + \\
- & - & - & - & - & + & + & + & + & - & - & + & + & - & + & + & - & + & - & - & + & + \\
- & + & + & - & + & - & - & + & + & - & + & + & + & - & - & + & + & + & + & - & - & + \\
+ & + & - & + & + & - & - & - & - & + & + & - & + & - & + & - & + & - & + & + & + & + \\
- & - & + & - & - & + & + & - & - & + & + & + & - & + & + & - & + & + & + & + & + & - \\
- & - & + & - & + & + & - & + & - & + & + & + & + & + & + & - & - & + & - & - & - & - \\
- & - & + & + & - & - & + & - & + & + & - & + & + & + & + & - & - & - & + & + & + & + \\
+ & - & + & - & + & - & + & + & + & - & + & - & + & + & + & - & + & + & + & + & + & + \\
+ & + & + & + & - & - & + & - & + & + & - & + & + & + & + & - & - & - & + & + & - & - \\
+ & + & + & - & + & + & + & + & - & - & - & - & + & - & + & + & + & + & - & - & + & - \\
+ & - & + & + & - & + & + & - & - & + & + & + & - & + & + & - & - & + & + & - & + & - \\
- & + & + & - & - & - & - & + & - & - & - & + & - & - & + & + & + & - & + & + & + & + \\
+ & + & + & - & + & + & - & + & + & - & + & + & - & - & - & + & + & - & + & - & - & + \\
+ & + & - & + & + & - & + & + & - & - & - & - & + & - & + & + & + & - & + & - & + & + \\
- & + & + & + & + & - & + & + & + & - & + & + & + & - & + & - & - & - & - & + & - & + \\
+ & - & + & - & + & + & - & - & + & + & - & - & + & + & + & + & + & + & - & - & + & - \\
+ & - & - & - & - & + & + & + & + & + & + & + & - & - & - & + & + & + & + & - & + & - \\
+ & - & - & + & + & + & - & - & + & + & + & - & + & + & + & - & + & - & - & + & + & - \\
\end{array}\right)\)}
\end{minipage}
\hfill
\begin{minipage}[t]{0.485\textwidth}
\centering
{\footnotesize $H_{2}=$}\par\vspace{-0.2mm}
{\fontsize{5.8}{5.9}\selectfont
\setlength{\arraycolsep}{0.55pt}\renewcommand{\arraystretch}{0.68}
\(\left(\begin{array}{*{22}{c}}
- & + & - & - & + & + & + & - & - & + & + & + & - & + & + & + & - & + & - & - & + & + \\
- & - & - & - & - & + & - & - & - & - & - & - & - & + & - & + & - & + & - & - & - & + \\
+ & - & + & - & + & + & - & + & - & - & - & + & + & - & + & - & - & + & - & - & + & + \\
- & - & - & - & - & + & + & + & + & + & - & + & + & - & + & - & + & + & + & + & - & + \\
- & - & - & - & + & + & - & - & - & + & + & - & + & + & + & - & - & + & + & + & + & - \\
+ & + & + & + & + & + & - & + & - & + & - & + & + & + & - & - & + & - & - & - & + & - \\
+ & - & - & - & - & - & + & + & + & + & - & - & + & + & - & + & + & - & + & - & - & + \\
+ & - & + & + & - & + & - & - & + & + & - & + & + & + & - & - & + & + & + & + & - & - \\
+ & + & + & - & + & + & - & - & - & - & + & + & - & + & - & + & - & + & - & + & + & + \\
- & - & - & + & - & - & + & + & - & - & + & + & + & - & + & + & - & + & + & + & + & + \\
- & - & - & + & - & + & + & - & + & - & + & + & + & + & + & + & - & - & + & - & - & - \\
+ & - & - & + & + & - & - & + & - & + & + & - & + & + & + & + & - & - & - & + & + & + \\
+ & + & - & + & - & + & - & + & + & + & - & + & - & + & + & + & - & + & + & + & + & + \\
- & + & + & + & + & - & - & + & - & + & + & - & + & + & + & + & - & - & - & + & + & - \\
- & + & + & + & - & + & + & + & + & - & - & - & - & + & - & + & + & + & + & - & - & + \\
- & + & - & + & + & - & + & + & - & - & + & + & + & - & + & + & - & - & + & + & - & + \\
+ & - & + & + & - & - & - & - & + & - & - & - & + & - & - & + & + & + & - & + & + & + \\
+ & + & + & + & - & + & + & - & + & + & - & + & + & - & - & - & + & + & - & + & - & - \\
+ & + & + & - & + & + & - & + & + & - & - & - & - & + & - & + & + & + & - & + & - & + \\
+ & - & + & + & + & + & - & + & + & + & - & + & + & + & - & + & - & - & - & - & + & - \\
- & + & - & + & - & + & + & - & - & + & + & - & - & + & + & + & + & + & + & - & - & + \\
- & + & - & - & - & - & + & + & + & + & + & + & + & - & - & - & + & + & + & + & - & + \\
\end{array}\right)\)}
\end{minipage}
\par\vspace{0.7mm}
\begin{minipage}[t]{0.485\textwidth}
\centering
{\footnotesize $H_{3}=$}\par\vspace{-0.2mm}
{\fontsize{5.8}{5.9}\selectfont
\setlength{\arraycolsep}{0.55pt}\renewcommand{\arraystretch}{0.68}
\(\left(\begin{array}{*{22}{c}}
+ & - & + & - & - & - & - & + & + & + & + & + & + & + & - & - & - & + & + & + & + & - \\
+ & - & + & - & - & + & + & + & - & - & + & + & + & - & + & + & + & - & + & - & - & + \\
+ & - & - & - & - & - & + & - & - & - & - & - & - & - & + & - & + & - & + & - & - & - \\
+ & + & - & + & - & + & + & - & + & - & - & - & + & + & - & + & - & - & + & - & - & + \\
+ & - & - & - & - & - & + & + & + & + & + & - & + & + & - & + & - & + & + & + & + & - \\
- & - & - & - & - & + & + & - & - & - & + & + & - & + & + & + & - & - & + & + & + & + \\
- & + & + & + & + & + & + & - & + & - & + & - & + & + & + & - & - & + & - & - & - & + \\
+ & + & - & - & - & - & - & + & + & + & + & - & - & + & + & - & + & + & - & + & - & - \\
- & + & - & + & + & - & + & - & - & + & + & - & + & + & + & - & - & + & + & + & + & - \\
+ & + & + & + & - & + & + & - & - & - & - & + & + & - & + & - & + & - & + & - & + & + \\
+ & - & - & - & + & - & - & + & + & - & - & + & + & + & - & + & + & - & + & + & + & + \\
- & - & - & - & + & - & + & + & - & + & - & + & + & + & + & + & + & - & - & + & - & - \\
+ & + & - & - & + & + & - & - & + & - & + & + & - & + & + & + & + & - & - & - & + & + \\
+ & + & + & - & + & - & + & - & + & + & + & - & + & - & + & + & + & - & + & + & + & + \\
- & - & + & + & + & + & - & - & + & - & + & + & - & + & + & + & + & - & - & - & + & + \\
+ & - & + & + & + & - & + & + & + & + & - & - & - & - & + & - & + & + & + & + & - & - \\
+ & - & + & - & + & + & - & + & + & - & - & + & + & + & - & + & + & - & - & + & + & - \\
+ & + & - & + & + & - & - & - & - & + & - & - & - & + & - & - & + & + & + & - & + & + \\
- & + & + & + & + & - & + & + & - & + & + & - & + & + & - & - & - & + & + & - & + & - \\
+ & + & + & + & - & + & + & - & + & + & - & - & - & - & + & - & + & + & + & - & + & - \\
- & + & - & + & + & + & + & - & + & + & + & - & + & + & + & - & + & - & - & - & - & + \\
+ & - & + & - & + & - & + & + & - & - & + & + & - & - & + & + & + & + & + & + & - & - \\
\end{array}\right)\)}
\end{minipage}
\hfill
\begin{minipage}[t]{0.485\textwidth}
\centering
{\footnotesize $H_{4}=$}\par\vspace{-0.2mm}
{\fontsize{5.8}{5.9}\selectfont
\setlength{\arraycolsep}{0.55pt}\renewcommand{\arraystretch}{0.68}
\(\left(\begin{array}{*{22}{c}}
- & + & - & + & - & + & - & + & + & - & - & + & + & - & - & + & + & + & + & + & + & - \\
- & + & - & + & - & - & - & - & + & + & + & + & + & + & + & - & - & - & + & + & + & + \\
+ & + & - & + & - & - & + & + & + & - & - & + & + & + & - & + & + & + & - & + & - & - \\
- & + & - & - & - & - & - & + & - & - & - & - & - & - & - & + & - & + & - & + & - & - \\
+ & + & + & - & + & - & + & + & - & + & - & - & - & + & + & - & + & - & - & + & - & - \\
- & + & - & - & - & - & - & + & + & + & + & + & - & + & + & - & + & - & + & + & + & + \\
+ & - & - & - & - & - & + & + & - & - & - & + & + & - & + & + & + & - & - & + & + & + \\
+ & - & + & + & + & + & + & + & - & + & - & + & - & + & + & + & - & - & + & - & - & - \\
- & + & + & - & - & - & - & - & + & + & + & + & - & - & + & + & - & + & + & - & + & - \\
- & - & + & - & + & + & - & + & - & - & + & + & - & + & + & + & - & - & + & + & + & + \\
+ & + & + & + & + & - & + & + & - & - & - & - & + & + & - & + & - & + & - & + & - & + \\
+ & + & - & - & - & + & - & - & + & + & - & - & + & + & + & - & + & + & - & + & + & + \\
- & - & - & - & - & + & - & + & + & - & + & - & + & + & + & + & + & + & - & - & + & - \\
+ & + & + & - & - & + & + & - & - & + & - & + & + & - & + & + & + & + & - & - & - & + \\
+ & + & + & + & - & + & - & + & - & + & + & + & - & + & - & + & + & + & - & + & + & + \\
+ & - & - & + & + & + & + & - & - & + & - & + & + & - & + & + & + & + & - & - & - & + \\
- & + & - & + & + & + & - & + & + & + & + & - & - & - & - & + & - & + & + & + & + & - \\
- & + & - & + & - & + & + & - & + & + & - & - & + & + & + & - & + & + & - & - & + & + \\
+ & + & + & - & + & + & - & - & - & - & + & - & - & - & + & - & - & + & + & + & - & + \\
- & - & + & + & + & + & - & + & + & - & + & + & - & + & + & - & - & - & + & + & - & + \\
- & + & + & + & + & - & + & + & - & + & + & - & - & - & - & + & - & + & + & + & - & + \\
+ & - & + & - & + & + & + & + & - & + & + & + & - & + & + & + & - & + & - & - & - & - \\
\end{array}\right)\)}
\end{minipage}
\par\vspace{0.7mm}
\begin{minipage}[t]{0.485\textwidth}
\centering
{\footnotesize $H_{5}=$}\par\vspace{-0.2mm}
{\fontsize{5.8}{5.9}\selectfont
\setlength{\arraycolsep}{0.55pt}\renewcommand{\arraystretch}{0.68}
\(\left(\begin{array}{*{22}{c}}
- & + & - & + & - & + & + & + & + & - & + & + & + & - & + & + & + & - & + & - & - & - \\
- & - & + & - & + & - & + & - & + & + & - & - & + & + & - & - & + & + & + & + & + & + \\
+ & - & + & - & + & - & - & - & - & + & + & + & + & + & + & + & - & - & - & + & + & + \\
- & + & + & - & + & - & - & + & + & + & - & - & + & + & + & - & + & + & + & - & + & - \\
- & - & + & - & - & - & - & - & + & - & - & - & - & - & - & - & + & - & + & - & + & - \\
- & + & + & + & - & + & - & + & + & - & + & - & - & - & + & + & - & + & - & - & + & - \\
+ & - & + & - & - & - & - & - & + & + & + & + & + & - & + & + & - & + & - & + & + & + \\
+ & + & - & - & - & - & - & + & + & - & - & - & + & + & - & + & + & + & - & - & + & + \\
- & + & - & + & + & + & + & + & + & - & + & - & + & - & + & + & + & - & - & + & - & - \\
- & - & + & + & - & - & - & - & - & + & + & + & + & - & - & + & + & - & + & + & - & + \\
+ & - & - & + & - & + & + & - & + & - & - & + & + & - & + & + & + & - & - & + & + & + \\
+ & + & + & + & + & + & - & + & + & - & - & - & - & + & + & - & + & - & + & - & + & - \\
+ & + & + & - & - & - & + & - & - & + & + & - & - & + & + & + & - & + & + & - & + & + \\
- & - & - & - & - & - & + & - & + & + & - & + & - & + & + & + & + & + & + & - & - & + \\
+ & + & + & + & - & - & + & + & - & - & + & - & + & + & - & + & + & + & + & - & - & - \\
+ & + & + & + & + & - & + & - & + & - & + & + & + & - & + & - & + & + & + & - & + & + \\
+ & + & - & - & + & + & + & + & - & - & + & - & + & + & - & + & + & + & + & - & - & - \\
- & - & + & - & + & + & + & - & + & + & + & + & - & - & - & - & + & - & + & + & + & + \\
+ & - & + & - & + & - & + & + & - & + & + & - & - & + & + & + & - & + & + & - & - & + \\
+ & + & + & + & - & + & + & - & - & - & - & + & - & - & - & + & - & - & + & + & + & - \\
+ & - & - & + & + & + & + & - & + & + & - & + & + & - & + & + & - & - & - & + & + & - \\
+ & - & + & + & + & + & - & + & + & - & + & + & - & - & - & - & + & - & + & + & + & - \\
\end{array}\right)\)}
\end{minipage}
\hfill
\begin{minipage}[t]{0.485\textwidth}
\centering
{\footnotesize $H_{6}=$}\par\vspace{-0.2mm}
{\fontsize{5.8}{5.9}\selectfont
\setlength{\arraycolsep}{0.55pt}\renewcommand{\arraystretch}{0.68}
\(\left(\begin{array}{*{22}{c}}
- & + & - & + & + & + & + & - & + & + & - & + & + & - & - & - & - & + & - & + & + & + \\
- & - & + & - & + & - & + & + & + & + & - & + & + & + & - & + & + & + & - & + & - & - \\
+ & - & - & + & - & + & - & + & - & + & + & - & - & + & + & - & - & + & + & + & + & + \\
+ & + & - & + & - & + & - & - & - & - & + & + & + & + & + & + & + & - & - & - & + & + \\
- & - & + & + & - & + & - & - & + & + & + & - & - & + & + & + & - & + & + & + & - & + \\
- & - & - & + & - & - & - & - & - & + & - & - & - & - & - & - & - & + & - & + & - & + \\
- & - & + & + & + & - & + & - & + & + & - & + & - & - & - & + & + & - & + & - & - & + \\
+ & + & - & + & - & - & - & - & - & + & + & + & + & + & - & + & + & - & + & - & + & + \\
+ & + & + & - & - & - & - & - & + & + & - & - & - & + & + & - & + & + & + & - & - & + \\
- & - & + & - & + & + & + & + & + & + & - & + & - & + & - & + & + & + & - & - & + & - \\
+ & - & - & + & + & - & - & - & - & - & + & + & + & + & - & - & + & + & - & + & + & - \\
+ & + & - & - & + & - & + & + & - & + & - & - & + & + & - & + & + & + & - & - & + & + \\
- & + & + & + & + & + & + & - & + & + & - & - & - & - & + & + & - & + & - & + & - & + \\
+ & + & + & + & - & - & - & + & - & - & + & + & - & - & + & + & + & - & + & + & - & + \\
+ & - & - & - & - & - & - & + & - & + & + & - & + & - & + & + & + & + & + & + & - & - \\
- & + & + & + & + & - & - & + & + & - & - & + & - & + & + & - & + & + & + & + & - & - \\
+ & + & + & + & + & + & - & + & - & + & - & + & + & + & - & + & - & + & + & + & - & + \\
- & + & + & - & - & + & + & + & + & - & - & + & - & + & + & - & + & + & + & + & - & - \\
+ & - & - & + & - & + & + & + & - & + & + & + & + & - & - & - & - & + & - & + & + & + \\
+ & + & - & + & - & + & - & + & + & - & + & + & - & - & + & + & + & - & + & + & - & - \\
- & + & + & + & + & - & + & + & - & - & - & - & + & - & - & - & + & - & - & + & + & + \\
- & + & - & - & + & + & + & + & - & + & + & - & + & + & - & + & + & - & - & - & + & + \\
\end{array}\right)\)}
\end{minipage}
\par\vspace{0.7mm}
\begin{minipage}[t]{0.485\textwidth}
\centering
{\footnotesize $H_{7}=$}\par\vspace{-0.2mm}
{\fontsize{5.8}{5.9}\selectfont
\setlength{\arraycolsep}{0.55pt}\renewcommand{\arraystretch}{0.68}
\(\left(\begin{array}{*{22}{c}}
+ & - & + & - & - & + & + & + & + & - & + & + & - & + & + & - & + & + & - & - & - & + \\
+ & - & + & - & + & + & + & + & - & + & + & - & + & + & - & - & - & - & + & - & + & + \\
- & - & - & + & - & + & - & + & + & + & + & - & + & + & + & - & + & + & + & - & + & - \\
+ & + & - & - & + & - & + & - & + & - & + & + & - & - & + & + & - & - & + & + & + & + \\
+ & + & + & - & + & - & + & - & - & - & - & + & + & + & + & + & + & + & - & - & - & + \\
+ & - & - & + & + & - & + & - & - & + & + & + & - & - & + & + & + & - & + & + & + & - \\
+ & - & - & - & + & - & - & - & - & - & + & - & - & - & - & - & - & - & + & - & + & - \\
+ & - & - & + & + & + & - & + & - & + & + & - & + & - & - & - & + & + & - & + & - & - \\
+ & + & + & - & + & - & - & - & - & - & + & + & + & + & + & - & + & + & - & + & - & + \\
+ & + & + & + & - & - & - & - & - & + & + & - & - & - & + & + & - & + & + & + & - & - \\
- & - & - & + & - & + & + & + & + & + & + & - & + & - & + & - & + & + & + & - & - & + \\
- & + & - & - & + & + & - & - & - & - & - & + & + & + & + & - & - & + & + & - & + & + \\
+ & + & + & - & - & + & - & + & + & - & + & - & - & + & + & - & + & + & + & - & - & + \\
+ & - & + & + & + & + & + & + & - & + & + & - & - & - & - & + & + & - & + & - & + & - \\
+ & + & + & + & + & - & - & - & + & - & - & + & + & - & - & + & + & + & - & + & + & - \\
- & + & - & - & - & - & - & - & + & - & + & + & - & + & - & + & + & + & + & + & + & - \\
- & - & + & + & + & + & - & - & + & + & - & - & + & - & + & + & - & + & + & + & + & - \\
+ & + & + & + & + & + & + & - & + & - & + & - & + & + & + & - & + & - & + & + & + & - \\
- & - & + & + & - & - & + & + & + & + & - & - & + & - & + & + & - & + & + & + & + & - \\
+ & + & - & - & + & - & + & + & + & - & + & + & + & + & - & - & - & - & + & - & + & + \\
- & + & + & - & + & - & + & - & + & + & - & + & + & - & - & + & + & + & - & + & + & - \\
+ & - & + & + & + & + & - & + & + & - & - & - & - & + & - & - & - & + & - & - & + & + \\
\end{array}\right)\)}
\end{minipage}
\hfill
\begin{minipage}[t]{0.485\textwidth}
\centering
{\footnotesize $H_{8}=$}\par\vspace{-0.2mm}
{\fontsize{5.8}{5.9}\selectfont
\setlength{\arraycolsep}{0.55pt}\renewcommand{\arraystretch}{0.68}
\(\left(\begin{array}{*{22}{c}}
+ & + & - & + & + & + & + & - & + & + & - & - & - & - & + & - & - & - & + & - & - & + \\
+ & + & - & + & - & - & + & + & + & + & - & + & + & - & + & + & - & + & + & - & - & - \\
+ & + & - & + & - & + & + & + & + & - & + & + & - & + & + & - & - & - & - & + & - & + \\
- & - & - & - & + & - & + & - & + & + & + & + & - & + & + & + & - & + & + & + & - & + \\
+ & + & + & - & - & + & - & + & - & + & - & + & + & - & - & + & + & - & - & + & + & + \\
+ & + & + & + & - & + & - & + & - & - & - & - & + & + & + & + & + & + & + & - & - & - \\
- & + & - & - & + & + & - & + & - & - & + & + & + & - & - & + & + & + & - & + & + & + \\
- & + & - & - & - & + & - & - & - & - & - & + & - & - & - & - & - & - & - & + & - & + \\
- & + & - & - & + & + & + & - & + & - & + & + & - & + & - & - & - & + & + & - & + & - \\
+ & + & + & + & - & + & - & - & - & - & - & + & + & + & + & + & - & + & + & - & + & - \\
- & + & + & + & + & - & - & - & - & - & + & + & - & - & - & + & + & - & + & + & + & - \\
+ & - & - & - & + & - & + & + & + & + & + & + & - & + & - & + & - & + & + & + & - & - \\
+ & - & + & - & - & + & + & - & - & - & - & - & + & + & + & + & - & - & + & + & - & + \\
+ & + & + & + & - & - & + & - & + & + & - & + & - & - & + & + & - & + & + & + & - & - \\
- & + & - & + & + & + & + & + & + & - & + & + & - & - & - & - & + & + & - & + & - & + \\
- & + & + & + & + & + & - & - & - & + & - & - & + & + & - & - & + & + & + & - & + & + \\
- & - & + & - & - & - & - & - & - & + & - & + & + & - & + & - & + & + & + & + & + & + \\
- & - & - & + & + & + & + & - & - & + & + & - & - & + & - & + & + & - & + & + & + & + \\
- & + & + & + & + & + & + & + & - & + & - & + & - & + & + & + & - & + & - & + & + & + \\
- & - & - & + & + & - & - & + & + & + & + & - & - & + & - & + & + & - & + & + & + & + \\
+ & + & + & - & - & + & - & + & + & + & - & + & + & + & + & - & - & - & - & + & - & + \\
- & - & + & + & - & + & - & + & - & + & + & - & + & + & - & - & + & + & + & - & + & + \\
\end{array}\right)\)}
\end{minipage}
\par\vspace{0.7mm}
\begin{minipage}[t]{0.485\textwidth}
\centering
{\footnotesize $H_{9}=$}\par\vspace{-0.2mm}
{\fontsize{5.8}{5.9}\selectfont
\setlength{\arraycolsep}{0.55pt}\renewcommand{\arraystretch}{0.68}
\(\left(\begin{array}{*{22}{c}}
+ & - & - & + & + & - & + & - & + & - & + & + & - & + & + & - & - & + & + & + & - & + \\
+ & + & + & - & + & + & + & + & - & + & + & - & - & - & - & + & - & - & - & + & - & - \\
- & + & + & - & + & - & - & + & + & + & + & - & + & + & - & + & + & - & + & + & - & - \\
+ & + & + & - & + & - & + & + & + & + & - & + & + & - & + & + & - & - & - & - & + & - \\
+ & - & - & - & - & + & - & + & - & + & + & + & + & - & + & + & + & - & + & + & + & - \\
+ & + & + & + & - & - & + & - & + & - & + & - & + & + & - & - & + & + & - & - & + & + \\
- & + & + & + & + & - & + & - & + & - & - & - & - & + & + & + & + & + & + & + & - & - \\
+ & - & + & - & - & + & + & - & + & - & - & + & + & + & - & - & + & + & + & - & + & + \\
+ & - & + & - & - & - & + & - & - & - & - & - & + & - & - & - & - & - & - & - & + & - \\
- & - & + & - & - & + & + & + & - & + & - & + & + & - & + & - & - & - & + & + & - & + \\
- & + & + & + & + & - & + & - & - & - & - & - & + & + & + & + & + & - & + & + & - & + \\
- & - & + & + & + & + & - & - & - & - & - & + & + & - & - & - & + & + & - & + & + & + \\
- & + & - & - & - & + & - & + & + & + & + & + & + & - & + & - & + & - & + & + & + & - \\
+ & + & - & + & - & - & + & + & - & - & - & - & - & + & + & + & + & - & - & + & + & - \\
- & + & + & + & + & - & - & + & - & + & + & - & + & - & - & + & + & - & + & + & + & - \\
+ & - & + & - & + & + & + & + & + & + & - & + & + & - & - & - & - & + & + & - & + & - \\
+ & - & + & + & + & + & + & - & - & - & + & - & - & + & + & - & - & + & + & + & - & + \\
+ & - & - & + & - & - & - & - & - & - & + & - & + & + & - & + & - & + & + & + & + & + \\
+ & - & - & - & + & + & + & + & - & - & + & + & - & - & + & - & + & + & - & + & + & + \\
+ & - & + & + & + & + & + & + & + & - & + & - & + & - & + & + & + & - & + & - & + & + \\
+ & - & - & - & + & + & - & - & + & + & + & + & - & - & + & - & + & + & - & + & + & + \\
+ & + & + & + & - & - & + & - & + & + & + & - & + & + & + & + & - & - & - & - & + & - \\
\end{array}\right)\)}
\end{minipage}
\hfill
\begin{minipage}[t]{0.485\textwidth}
\centering
{\footnotesize $H_{10}=$}\par\vspace{-0.2mm}
{\fontsize{5.8}{5.9}\selectfont
\setlength{\arraycolsep}{0.55pt}\renewcommand{\arraystretch}{0.68}
\(\left(\begin{array}{*{22}{c}}
- & + & + & + & + & - & - & + & - & + & + & + & - & + & + & + & + & - & - & - & - & + \\
+ & + & - & - & + & + & - & + & - & + & - & + & + & - & + & + & - & - & + & + & + & - \\
- & + & + & + & - & + & + & + & + & - & + & + & - & - & - & - & + & - & - & - & + & - \\
- & - & + & + & - & + & - & - & + & + & + & + & - & + & + & - & + & + & - & + & + & - \\
- & + & + & + & - & + & - & + & + & + & + & - & + & + & - & + & + & - & - & - & - & + \\
- & + & - & - & - & - & + & - & + & - & + & + & + & + & - & + & + & + & - & + & + & + \\
+ & + & + & + & + & - & - & + & - & + & - & + & - & + & + & - & - & + & + & - & - & + \\
- & - & + & + & + & + & - & + & - & + & - & - & - & - & + & + & + & + & + & + & + & - \\
+ & + & - & + & - & - & + & + & - & + & - & - & + & + & + & - & - & + & + & + & - & + \\
- & + & - & + & - & - & - & + & - & - & - & - & - & + & - & - & - & - & - & - & - & + \\
+ & - & - & + & - & - & + & + & + & - & + & - & + & + & - & + & - & - & - & + & + & - \\
+ & - & + & + & + & + & - & + & - & - & - & - & - & + & + & + & + & + & - & + & + & - \\
+ & - & - & + & + & + & + & - & - & - & - & - & + & + & - & - & - & + & + & - & + & + \\
- & - & + & - & - & - & + & - & + & + & + & + & + & + & - & + & - & + & - & + & + & + \\
- & + & + & - & + & - & - & + & + & - & - & - & - & - & + & + & + & + & - & - & + & + \\
- & - & + & + & + & + & - & - & + & - & + & + & - & + & - & - & + & + & - & + & + & + \\
- & + & - & + & - & + & + & + & + & + & + & - & + & + & - & - & - & - & + & + & - & + \\
+ & + & - & + & + & + & + & + & - & - & - & + & - & - & + & + & - & - & + & + & + & - \\
+ & + & - & - & + & - & - & - & - & - & - & + & - & + & + & - & + & - & + & + & + & + \\
+ & + & - & - & - & + & + & + & + & - & - & + & + & - & - & + & - & + & + & - & + & + \\
+ & + & - & + & + & + & + & + & + & + & - & + & - & + & - & + & + & + & - & + & - & + \\
+ & + & - & - & - & + & + & - & - & + & + & + & + & - & - & + & - & + & + & - & + & + \\
\end{array}\right)\)}
\end{minipage}
\par\vspace{0.7mm}
\begin{minipage}[t]{0.485\textwidth}
\centering
{\footnotesize $H_{11}=$}\par\vspace{-0.2mm}
{\fontsize{5.8}{5.9}\selectfont
\setlength{\arraycolsep}{0.55pt}\renewcommand{\arraystretch}{0.68}
\(\left(\begin{array}{*{22}{c}}
+ & + & + & - & - & - & + & + & - & - & + & + & + & + & - & - & + & - & + & + & - & + \\
+ & - & + & + & + & + & - & - & + & - & + & + & + & - & + & + & + & + & - & - & - & - \\
- & + & + & - & - & + & + & - & + & - & + & - & + & + & - & + & + & - & - & + & + & + \\
- & - & + & + & + & - & + & + & + & + & - & + & + & - & - & - & - & + & - & - & - & + \\
- & - & - & + & + & - & + & - & - & + & + & + & + & - & + & + & - & + & + & - & + & + \\
+ & - & + & + & + & - & + & - & + & + & + & + & - & + & + & - & + & + & - & - & - & - \\
+ & - & + & - & - & - & - & + & - & + & - & + & + & + & + & - & + & + & + & - & + & + \\
+ & + & + & + & + & + & - & - & + & - & + & - & + & - & + & + & - & - & + & + & - & - \\
- & - & - & + & + & + & + & - & + & - & + & - & - & - & - & + & + & + & + & + & + & + \\
+ & + & + & - & + & - & - & + & + & - & + & - & - & + & + & + & - & - & + & + & + & - \\
+ & - & + & - & + & - & - & - & + & - & - & - & - & - & + & - & - & - & - & - & - & - \\
- & + & - & - & + & - & - & + & + & + & - & + & - & + & + & - & + & - & - & - & + & + \\
- & + & - & + & + & + & + & - & + & - & - & - & - & - & + & + & + & + & + & - & + & + \\
+ & + & - & - & + & + & + & + & - & - & - & - & - & + & + & - & - & - & + & + & - & + \\
+ & - & - & + & - & - & - & + & - & + & + & + & + & + & + & - & + & - & + & - & + & + \\
+ & - & + & + & - & + & - & - & + & + & - & - & - & - & - & + & + & + & + & - & - & + \\
+ & - & - & + & + & + & + & - & - & + & - & + & + & - & + & - & - & + & + & - & + & + \\
+ & - & + & - & + & - & + & + & + & + & + & + & - & + & + & - & - & - & - & + & + & - \\
- & + & + & - & + & + & + & + & + & - & - & - & + & - & - & + & + & - & - & + & + & + \\
+ & + & + & - & - & + & - & - & - & - & - & - & + & - & + & + & - & + & - & + & + & + \\
+ & + & + & - & - & - & + & + & + & + & - & - & + & + & - & - & + & - & + & + & - & + \\
+ & + & + & - & + & + & + & + & + & + & + & - & + & - & + & - & + & + & + & - & + & - \\
\end{array}\right)\)}
\end{minipage}
\hfill
\begin{minipage}[t]{0.485\textwidth}
\centering
{\footnotesize $H_{12}=$}\par\vspace{-0.2mm}
{\fontsize{5.8}{5.9}\selectfont
\setlength{\arraycolsep}{0.55pt}\renewcommand{\arraystretch}{0.68}
\(\left(\begin{array}{*{22}{c}}
- & + & + & + & - & + & + & + & + & + & + & + & - & + & - & + & - & + & + & + & - & + \\
+ & + & + & + & - & - & - & + & + & - & - & + & + & + & + & - & - & + & - & + & + & - \\
- & + & - & + & + & + & + & - & - & + & - & + & + & + & - & + & + & + & + & - & - & - \\
+ & - & + & + & - & - & + & + & - & + & - & + & - & + & + & - & + & + & - & - & + & + \\
+ & - & - & + & + & + & - & + & + & + & + & - & + & + & - & - & - & - & + & - & - & - \\
+ & - & - & - & + & + & - & + & - & - & + & + & + & + & - & + & + & - & + & + & - & + \\
- & + & - & + & + & + & - & + & - & + & + & + & + & - & + & + & - & + & + & - & - & - \\
+ & + & - & + & - & - & - & - & + & - & + & - & + & + & + & + & - & + & + & + & - & + \\
- & + & + & + & + & + & + & - & - & + & - & + & - & + & - & + & + & - & - & + & + & - \\
+ & - & - & - & + & + & + & + & - & + & - & + & - & - & - & - & + & + & + & + & + & + \\
- & + & + & + & - & + & - & - & + & + & - & + & - & - & + & + & + & - & - & + & + & + \\
- & + & - & + & - & + & - & - & - & + & - & - & - & - & - & + & - & - & - & - & - & - \\
+ & - & + & - & - & + & - & - & + & + & + & - & + & - & + & + & - & + & - & - & - & + \\
+ & - & + & - & + & + & + & + & - & + & - & - & - & - & - & + & + & + & + & + & - & + \\
+ & + & + & - & - & + & + & + & + & - & - & - & - & - & + & + & - & - & - & + & + & - \\
+ & + & - & - & + & - & - & - & + & - & + & + & + & + & + & + & - & + & - & + & - & + \\
+ & + & - & + & + & - & + & - & - & + & + & - & - & - & - & - & + & + & + & + & - & - \\
+ & + & - & - & + & + & + & + & - & - & + & - & + & + & - & + & - & - & + & + & - & + \\
- & + & - & + & - & + & - & + & + & + & + & + & + & - & + & + & - & - & - & - & + & + \\
+ & - & + & + & - & + & + & + & + & + & - & - & - & + & - & - & + & + & - & - & + & + \\
+ & + & + & + & - & - & + & - & - & - & - & - & - & + & - & + & + & - & + & - & + & + \\
+ & + & + & + & - & - & - & + & + & + & + & - & - & + & + & - & - & + & - & + & + & - \\
\end{array}\right)\)}
\end{minipage}
\clearpage
\begin{minipage}[t]{0.485\textwidth}
\centering
{\footnotesize $H_{13}=$}\par\vspace{-0.2mm}
{\fontsize{5.8}{5.9}\selectfont
\setlength{\arraycolsep}{0.55pt}\renewcommand{\arraystretch}{0.68}
\(\left(\begin{array}{*{22}{c}}
- & + & + & + & + & - & - & - & + & + & + & + & - & - & + & + & - & - & + & - & + & + \\
+ & - & + & + & + & - & + & + & + & + & + & + & + & - & + & - & + & - & + & + & + & - \\
- & + & + & + & + & - & - & - & + & + & - & - & + & + & + & + & - & - & + & - & + & + \\
- & - & + & - & + & + & + & + & - & - & + & - & + & + & + & - & + & + & + & + & - & - \\
+ & + & - & + & + & - & - & + & + & - & + & - & + & - & + & + & - & + & + & - & - & + \\
- & + & - & - & + & + & + & - & + & + & + & + & - & + & + & - & - & - & - & + & - & - \\
+ & + & - & - & - & + & + & - & + & - & - & + & + & + & + & - & + & + & - & + & + & - \\
- & - & + & - & + & + & + & - & + & - & + & + & + & + & - & + & + & - & + & + & - & - \\
+ & + & + & - & + & - & - & - & - & + & - & + & - & + & + & + & + & - & + & + & + & - \\
- & - & + & + & + & + & + & + & - & - & + & - & + & - & + & - & + & + & - & - & + & + \\
+ & + & - & - & - & + & + & + & + & - & + & - & + & - & - & - & - & + & + & + & + & + \\
+ & - & + & + & + & - & + & - & - & + & + & - & + & - & - & + & + & + & - & - & + & + \\
- & - & + & - & + & - & + & - & - & - & + & - & - & - & - & - & + & - & - & - & - & - \\
+ & + & - & + & - & - & + & - & - & + & + & + & - & + & - & + & + & - & + & - & - & - \\
+ & + & - & + & - & + & + & + & + & - & + & - & - & - & - & - & + & + & + & + & + & - \\
- & + & + & + & - & - & + & + & + & + & - & - & - & - & - & + & + & - & - & - & + & + \\
+ & + & + & - & - & + & - & - & - & + & - & + & + & + & + & + & + & - & + & - & + & - \\
- & + & + & - & + & + & - & + & - & - & + & + & - & - & - & - & - & + & + & + & + & - \\
+ & + & + & - & - & + & + & + & + & - & - & + & - & + & + & - & + & - & - & + & + & - \\
+ & - & + & - & + & - & + & - & + & + & + & + & + & + & - & + & + & - & - & - & - & + \\
+ & + & - & + & + & - & + & + & + & + & + & - & - & - & + & - & - & + & + & - & - & + \\
+ & + & + & + & + & - & - & + & - & - & - & - & - & - & + & - & + & + & - & + & - & + \\
\end{array}\right)\)}
\end{minipage}
\hfill
\begin{minipage}[t]{0.485\textwidth}
\centering
{\footnotesize $H_{14}=$}\par\vspace{-0.2mm}
{\fontsize{5.8}{5.9}\selectfont
\setlength{\arraycolsep}{0.55pt}\renewcommand{\arraystretch}{0.68}
\(\left(\begin{array}{*{22}{c}}
+ & + & + & + & + & + & - & - & + & - & - & - & - & - & - & + & - & + & + & - & + & - \\
+ & - & + & + & + & + & - & - & - & + & + & + & + & - & - & + & + & - & - & + & - & + \\
- & + & - & + & + & + & - & + & + & + & + & + & + & + & - & + & - & + & - & + & + & + \\
+ & - & + & + & + & + & - & - & - & + & + & - & - & + & + & + & + & - & - & + & - & + \\
- & - & - & + & - & + & + & + & + & - & - & + & - & + & + & + & - & + & + & + & + & - \\
+ & + & + & - & + & + & - & - & + & + & - & + & - & + & - & + & + & - & + & + & - & - \\
- & - & + & - & - & + & + & + & - & + & + & + & + & - & + & + & - & - & - & - & + & - \\
- & + & + & - & - & - & + & + & - & + & - & - & + & + & + & + & - & + & + & - & + & + \\
- & - & - & + & - & + & + & + & - & + & - & + & + & + & + & - & + & + & - & + & + & - \\
- & + & + & + & - & + & - & - & - & - & + & - & + & - & + & + & + & + & - & + & + & + \\
+ & - & - & + & + & + & + & + & + & - & - & + & - & + & - & + & - & + & + & - & - & + \\
+ & + & + & - & - & - & + & + & + & + & - & + & - & + & - & - & - & - & + & + & + & + \\
+ & + & - & + & + & + & - & + & - & - & + & + & - & + & - & - & + & + & + & - & - & + \\
- & - & - & + & - & + & - & + & - & - & - & + & - & - & - & - & - & + & - & - & - & - \\
- & + & + & - & + & - & - & + & - & - & + & + & + & - & + & - & + & + & - & + & - & - \\
- & + & + & - & + & - & + & + & + & + & - & + & - & - & - & - & - & + & + & + & + & + \\
+ & - & + & + & + & - & - & + & + & + & + & - & - & - & - & - & + & + & - & - & - & + \\
- & + & + & + & - & - & + & - & - & - & + & - & + & + & + & + & + & + & - & + & - & + \\
- & - & + & + & - & + & + & - & + & - & - & + & + & - & - & - & - & - & + & + & + & + \\
- & + & + & + & - & - & + & + & + & + & - & - & + & - & + & + & - & + & - & - & + & + \\
+ & + & - & + & - & + & - & + & - & + & + & + & + & + & + & - & + & + & - & - & - & - \\
+ & + & + & - & + & + & - & + & + & + & + & + & - & - & - & + & - & - & + & + & - & - \\
\end{array}\right)\)}
\end{minipage}
\par\vspace{0.7mm}
\begin{minipage}[t]{0.485\textwidth}
\centering
{\footnotesize $H_{15}=$}\par\vspace{-0.2mm}
{\fontsize{5.8}{5.9}\selectfont
\setlength{\arraycolsep}{0.55pt}\renewcommand{\arraystretch}{0.68}
\(\left(\begin{array}{*{22}{c}}
- & + & + & + & - & + & + & - & + & + & + & + & + & - & - & - & + & - & - & + & + & - \\
- & + & + & + & + & + & + & - & - & + & - & - & - & - & - & - & + & - & + & + & - & + \\
+ & + & - & + & + & + & + & - & - & - & + & + & + & + & - & - & + & + & - & - & + & - \\
+ & - & + & - & + & + & + & - & + & + & + & + & + & + & + & - & + & - & + & - & + & + \\
+ & + & - & + & + & + & + & - & - & - & + & + & - & - & + & + & + & + & - & - & + & - \\
- & - & - & - & + & - & + & + & + & + & - & - & + & - & + & + & + & - & + & + & + & + \\
- & + & + & + & - & + & + & - & - & + & + & - & + & - & + & - & + & + & - & + & + & - \\
- & - & - & + & - & - & + & + & + & - & + & + & + & + & - & + & + & - & - & - & - & + \\
+ & - & + & + & - & - & - & + & + & - & + & - & - & + & + & + & + & - & + & + & - & + \\
- & - & - & - & + & - & + & + & + & - & + & - & + & + & + & + & - & + & + & - & + & + \\
+ & - & + & + & + & - & + & - & - & - & - & + & - & + & - & + & + & + & + & - & + & + \\
+ & + & - & - & + & + & + & + & + & + & - & - & + & - & + & - & + & - & + & + & - & - \\
+ & + & + & + & - & - & - & + & + & + & + & - & + & - & + & - & - & - & - & + & + & + \\
+ & + & + & - & + & + & + & - & + & - & - & + & + & - & + & - & - & + & + & + & - & - \\
- & - & - & - & + & - & + & - & + & - & - & - & + & - & - & - & - & - & + & - & - & - \\
- & - & + & + & - & + & - & - & + & - & - & + & + & + & - & + & - & + & + & - & + & - \\
+ & - & + & + & - & + & - & + & + & + & + & - & + & - & - & - & - & - & + & + & + & + \\
+ & + & - & + & + & + & - & - & + & + & + & + & - & - & - & - & - & + & + & - & - & - \\
+ & - & + & + & + & - & - & + & - & - & - & + & - & + & + & + & + & + & + & - & + & - \\
+ & - & - & + & + & - & + & + & - & + & - & - & + & + & - & - & - & - & - & + & + & + \\
+ & - & + & + & + & - & - & + & + & + & + & - & - & + & - & + & + & - & + & - & - & + \\
- & + & + & - & + & - & + & - & + & - & + & + & + & + & + & + & - & + & + & - & - & - \\
\end{array}\right)\)}
\end{minipage}
\hfill
\begin{minipage}[t]{0.485\textwidth}
\centering
{\footnotesize $H_{16}=$}\par\vspace{-0.2mm}
{\fontsize{5.8}{5.9}\selectfont
\setlength{\arraycolsep}{0.55pt}\renewcommand{\arraystretch}{0.68}
\(\left(\begin{array}{*{22}{c}}
- & - & + & + & - & + & - & + & - & + & - & + & + & + & + & + & + & - & + & + & - & - \\
- & - & + & + & + & - & + & + & - & + & + & + & + & + & - & - & - & + & - & - & + & + \\
+ & - & + & + & + & + & + & + & - & - & + & - & - & - & - & - & - & + & - & + & + & - \\
- & + & + & - & + & + & + & + & - & - & - & + & + & + & + & - & - & + & + & - & - & + \\
+ & + & - & + & - & + & + & + & - & + & + & + & + & + & + & + & - & + & - & + & - & + \\
- & + & + & - & + & + & + & + & - & - & - & + & + & - & - & + & + & + & + & - & - & + \\
+ & - & - & - & - & + & - & + & + & + & + & - & - & + & - & + & + & + & - & + & + & + \\
- & - & + & + & + & - & + & + & - & - & + & + & - & + & - & + & - & + & + & - & + & + \\
+ & - & - & - & + & - & - & + & + & + & - & + & + & + & + & - & + & + & - & - & - & - \\
+ & + & - & + & + & - & - & - & + & + & - & + & - & - & + & + & + & + & - & + & + & - \\
+ & - & - & - & - & + & - & + & + & + & - & + & - & + & + & + & + & - & + & + & - & + \\
+ & + & - & + & + & + & - & + & - & - & - & - & + & - & + & - & + & + & + & + & - & + \\
- & + & + & - & - & + & + & + & + & + & + & - & - & + & - & + & - & + & - & + & + & - \\
+ & + & + & + & + & - & - & - & + & + & + & + & - & + & - & + & - & - & - & - & + & + \\
- & + & + & + & - & + & + & + & - & + & - & - & + & + & - & + & - & - & + & + & + & - \\
- & - & - & - & - & + & - & + & - & + & - & - & - & + & - & - & - & - & - & + & - & - \\
- & - & - & + & + & - & + & - & - & + & - & - & + & + & + & - & + & - & + & + & - & + \\
+ & + & - & + & + & - & + & - & + & + & + & + & - & + & - & - & - & - & - & + & + & + \\
- & + & + & - & + & + & + & - & - & + & + & + & + & - & - & - & - & - & + & + & - & - \\
- & + & - & + & + & + & - & - & + & - & - & - & + & - & + & + & + & + & + & + & - & + \\
+ & + & - & - & + & + & - & + & + & - & + & - & - & + & + & - & - & - & - & - & + & + \\
+ & + & - & + & + & + & - & - & + & + & + & + & - & - & + & - & + & + & - & + & - & - \\
\end{array}\right)\)}
\end{minipage}
\par\vspace{0.7mm}
\begin{minipage}[t]{0.485\textwidth}
\centering
{\footnotesize $H_{17}=$}\par\vspace{-0.2mm}
{\fontsize{5.8}{5.9}\selectfont
\setlength{\arraycolsep}{0.55pt}\renewcommand{\arraystretch}{0.68}
\(\left(\begin{array}{*{22}{c}}
- & + & + & - & + & + & + & - & - & + & + & + & + & - & - & + & - & + & + & - & + & - \\
- & - & - & + & + & - & + & - & + & - & + & - & + & + & + & + & + & + & - & + & + & - \\
+ & - & - & + & + & + & - & + & + & - & + & + & + & + & + & - & - & - & + & - & - & + \\
- & + & - & + & + & + & + & + & + & - & - & + & - & - & - & - & - & - & + & - & + & + \\
+ & - & + & + & - & + & + & + & + & - & - & - & + & + & + & + & - & - & + & + & - & - \\
+ & + & + & - & + & - & + & + & + & - & + & + & + & + & + & + & + & - & + & - & + & - \\
+ & - & + & + & - & + & + & + & + & - & - & - & + & + & - & - & + & + & + & + & - & - \\
+ & + & - & - & - & - & + & - & + & + & + & + & - & - & + & - & + & + & + & - & + & + \\
+ & - & - & + & + & + & - & + & + & - & - & + & + & - & + & - & + & - & + & + & - & + \\
- & + & - & - & - & + & - & - & + & + & + & - & + & + & + & + & - & + & + & - & - & - \\
- & + & + & - & + & + & - & - & - & + & + & - & + & - & - & + & + & + & + & - & + & + \\
+ & + & - & - & - & - & + & - & + & + & + & - & + & - & + & + & + & + & - & + & + & - \\
+ & + & + & - & + & + & + & - & + & - & - & - & - & + & - & + & - & + & + & + & + & - \\
- & - & + & + & - & - & + & + & + & + & + & + & - & - & + & - & + & - & + & - & + & + \\
+ & + & + & + & + & + & - & - & - & + & + & + & + & - & + & - & + & - & - & - & - & + \\
- & - & + & + & + & - & + & + & + & - & + & - & - & + & + & - & + & - & - & + & + & + \\
- & - & - & - & - & - & + & - & + & - & + & - & - & - & + & - & - & - & - & - & + & - \\
+ & - & - & - & + & + & - & + & - & - & + & - & - & + & + & + & - & + & - & + & + & - \\
+ & + & + & - & + & + & - & + & - & + & + & + & + & - & + & - & - & - & - & - & + & + \\
- & - & + & + & - & + & + & + & - & - & + & + & + & + & - & - & - & - & - & + & + & - \\
+ & - & + & - & + & + & + & - & - & + & - & - & - & + & - & + & + & + & + & + & + & - \\
+ & + & + & - & - & + & + & - & + & + & - & + & - & - & + & + & - & - & - & - & - & + \\
\end{array}\right)\)}
\end{minipage}
\hfill
\begin{minipage}[t]{0.485\textwidth}
\centering
{\footnotesize $H_{18}=$}\par\vspace{-0.2mm}
{\fontsize{5.8}{5.9}\selectfont
\setlength{\arraycolsep}{0.55pt}\renewcommand{\arraystretch}{0.68}
\(\left(\begin{array}{*{22}{c}}
+ & + & + & + & - & - & + & + & - & + & + & - & + & - & - & + & + & - & - & - & - & - \\
- & - & + & + & - & + & + & + & - & - & + & + & + & + & - & - & + & - & + & + & - & + \\
- & - & - & - & + & + & - & + & - & + & - & + & - & + & + & + & + & + & + & - & + & + \\
+ & + & - & - & + & + & + & - & + & + & - & + & + & + & + & + & - & - & - & + & - & - \\
+ & - & + & - & + & + & + & + & + & + & - & - & + & - & - & - & - & - & - & + & - & + \\
- & + & - & + & + & - & + & + & + & + & - & - & - & + & + & + & + & - & - & + & + & - \\
- & + & + & + & - & + & - & + & + & + & - & + & + & + & + & + & + & + & - & + & - & + \\
- & + & - & + & + & - & + & + & + & + & - & - & - & + & + & - & - & + & + & + & + & - \\
+ & + & + & - & - & - & - & + & - & + & + & + & + & - & - & + & - & + & + & + & - & + \\
+ & + & - & - & + & + & + & - & + & + & - & - & + & + & - & + & - & + & - & + & + & - \\
- & - & + & - & - & - & + & - & - & + & + & + & - & + & + & + & + & - & + & + & - & - \\
+ & - & + & + & - & + & + & - & - & - & + & + & - & + & - & - & + & + & + & + & - & + \\
- & + & + & - & - & - & - & + & - & + & + & + & - & + & - & + & + & + & + & - & + & + \\
- & + & + & + & - & + & + & + & - & + & - & - & - & - & + & - & + & - & + & + & + & + \\
+ & - & - & + & + & - & - & + & + & + & + & + & + & - & - & + & - & + & - & + & - & + \\
+ & + & + & + & + & + & + & - & - & - & + & + & + & + & - & + & - & + & - & - & - & - \\
+ & - & - & + & + & + & - & + & + & + & - & + & - & - & + & + & - & + & - & - & + & + \\
- & - & - & - & - & - & - & + & - & + & - & + & - & - & - & + & - & - & - & - & - & + \\
- & + & - & - & - & + & + & - & + & - & - & + & - & - & + & + & + & - & + & - & + & + \\
+ & + & + & + & - & + & + & - & + & - & + & + & + & + & - & + & - & - & - & - & - & + \\
- & - & - & + & + & - & + & + & + & - & - & + & + & + & + & - & - & - & - & - & + & + \\
- & + & - & + & - & + & + & + & - & - & + & - & - & - & + & - & + & + & + & + & + & + \\
\end{array}\right)\)}
\end{minipage}
\par\vspace{0.7mm}
\begin{minipage}[t]{0.485\textwidth}
\centering
{\footnotesize $H_{19}=$}\par\vspace{-0.2mm}
{\fontsize{5.8}{5.9}\selectfont
\setlength{\arraycolsep}{0.55pt}\renewcommand{\arraystretch}{0.68}
\(\left(\begin{array}{*{22}{c}}
+ & - & + & - & + & - & + & + & + & - & - & + & - & - & - & + & - & + & + & + & + & + \\
- & + & + & + & + & - & - & + & + & - & + & + & - & + & - & - & + & + & - & - & - & - \\
+ & - & - & + & + & - & + & + & + & - & - & + & + & + & + & - & - & + & - & + & + & - \\
+ & - & - & - & - & + & + & - & + & - & + & - & + & - & + & + & + & + & + & + & - & + \\
- & + & + & - & - & + & + & + & - & + & + & - & + & + & + & + & + & - & - & - & + & - \\
+ & + & - & + & - & + & + & + & + & + & + & - & - & + & - & - & - & - & - & - & + & - \\
- & - & + & - & + & + & - & + & + & + & + & - & - & - & + & + & + & + & - & - & + & + \\
+ & - & + & + & + & - & + & - & + & + & + & - & + & + & + & + & + & + & + & - & + & - \\
- & - & + & - & + & + & - & + & + & + & + & - & - & - & + & + & - & - & + & + & + & + \\
+ & + & + & + & - & - & - & - & + & - & + & + & + & + & - & - & + & - & + & + & + & - \\
- & + & + & - & - & + & + & + & - & + & + & - & - & + & + & - & + & - & + & - & + & + \\
- & - & - & + & - & - & - & + & - & - & + & + & + & - & + & + & + & + & - & + & + & - \\
+ & + & - & + & + & - & + & + & - & - & - & + & + & - & + & - & - & + & + & + & + & - \\
+ & - & + & + & - & - & - & - & + & - & + & + & + & - & + & - & + & + & + & + & - & + \\
+ & - & + & + & + & - & + & + & + & - & + & - & - & - & - & + & - & + & - & + & + & + \\
+ & + & - & - & + & + & - & - & + & + & + & + & + & + & - & - & + & - & + & - & + & - \\
- & + & + & + & + & + & + & + & - & - & - & + & + & + & + & - & + & - & + & - & - & - \\
+ & + & - & - & + & + & + & - & + & + & + & - & + & - & - & + & + & - & + & - & - & + \\
+ & - & - & - & - & - & - & - & + & - & + & - & + & - & - & - & + & - & - & - & - & - \\
+ & - & + & - & - & - & + & + & - & + & - & - & + & - & - & + & + & + & - & + & - & + \\
+ & + & + & + & + & - & + & + & - & + & - & + & + & + & + & - & + & - & - & - & - & - \\
+ & - & - & - & + & + & - & + & + & + & - & - & + & + & + & + & - & - & - & - & - & + \\
\end{array}\right)\)}
\end{minipage}
\hfill
\begin{minipage}[t]{0.485\textwidth}
\centering
{\footnotesize $H_{20}=$}\par\vspace{-0.2mm}
{\fontsize{5.8}{5.9}\selectfont
\setlength{\arraycolsep}{0.55pt}\renewcommand{\arraystretch}{0.68}
\(\left(\begin{array}{*{22}{c}}
+ & + & - & - & - & + & + & - & + & + & + & - & - & + & + & + & + & - & - & - & - & - \\
+ & + & - & + & - & + & - & + & + & + & - & - & + & - & - & - & + & - & + & + & + & + \\
- & - & + & + & + & + & - & - & + & + & - & + & + & - & + & - & - & + & + & - & - & - \\
- & + & - & - & + & + & - & + & + & + & - & - & + & + & + & + & - & - & + & - & + & + \\
+ & + & - & - & - & - & + & + & - & + & - & + & - & + & - & + & + & + & + & + & + & - \\
- & - & + & + & - & - & + & + & + & - & + & + & - & + & + & + & + & + & - & - & - & + \\
- & + & + & - & + & - & + & + & + & + & + & + & - & - & + & - & - & - & - & - & - & + \\
+ & - & - & + & - & + & + & - & + & + & + & + & - & - & - & + & + & + & + & - & - & + \\
- & + & - & + & + & + & - & + & - & + & + & + & - & + & + & + & + & + & + & + & - & + \\
+ & - & - & + & - & + & + & - & + & + & + & + & - & - & - & + & + & - & - & + & + & + \\
- & + & + & + & + & - & - & - & - & + & - & + & + & + & + & - & - & + & - & + & + & + \\
+ & - & + & + & - & - & + & + & + & - & + & + & - & - & + & + & - & + & - & + & - & + \\
- & - & - & - & + & - & - & - & + & - & - & + & + & + & - & + & + & + & + & - & + & + \\
- & + & + & - & + & + & - & + & + & - & - & - & + & + & - & + & - & - & + & + & + & + \\
+ & + & - & + & + & - & - & - & - & + & - & + & + & + & - & + & - & + & + & + & + & - \\
+ & + & - & + & + & + & - & + & + & + & - & + & - & - & - & - & + & - & + & - & + & + \\
- & + & + & - & - & + & + & - & - & + & + & + & + & + & + & - & - & + & - & + & - & + \\
- & - & + & + & + & + & + & + & + & - & - & - & + & + & + & + & - & + & - & + & - & - \\
+ & + & + & - & - & + & + & + & - & + & + & + & - & + & - & - & + & + & - & + & - & - \\
- & + & - & - & - & - & - & - & - & + & - & + & - & + & - & - & - & + & - & - & - & - \\
+ & + & - & + & - & - & - & + & + & - & + & - & - & + & - & - & + & + & + & - & + & - \\
- & + & + & + & + & + & - & + & + & - & + & - & + & + & + & + & - & + & - & - & - & - \\
\end{array}\right)\)}
\end{minipage}
\par\vspace{0.7mm}
\begin{minipage}[t]{0.485\textwidth}
\centering
{\footnotesize $H_{21}=$}\par\vspace{-0.2mm}
{\fontsize{5.8}{5.9}\selectfont
\setlength{\arraycolsep}{0.55pt}\renewcommand{\arraystretch}{0.68}
\(\left(\begin{array}{*{22}{c}}
- & - & + & + & + & + & + & - & + & + & - & + & - & + & + & + & + & - & + & - & - & - \\
- & + & + & - & - & - & + & + & - & + & + & + & - & - & + & + & + & + & - & - & - & - \\
+ & + & + & - & + & - & + & - & + & + & + & - & - & + & - & - & - & + & - & + & + & + \\
- & - & - & + & + & + & + & - & - & + & + & - & + & + & - & + & - & - & + & + & - & - \\
+ & - & + & - & - & + & + & - & + & + & + & - & - & + & + & + & + & - & - & + & - & + \\
- & + & + & - & - & - & - & + & + & - & + & - & + & - & + & - & + & + & + & + & + & + \\
+ & - & - & + & + & - & - & + & + & + & - & + & + & - & + & + & + & + & + & - & - & - \\
+ & - & + & + & - & + & - & + & + & + & + & + & + & - & - & + & - & - & - & - & - & - \\
+ & + & - & - & + & - & + & + & - & + & + & + & + & - & - & - & + & + & + & + & - & - \\
+ & - & + & - & + & + & + & - & + & - & + & + & + & - & + & + & + & + & + & + & + & - \\
+ & + & - & - & + & - & + & + & - & + & + & + & + & - & - & - & + & + & - & - & + & + \\
+ & - & + & + & + & + & - & - & - & - & + & - & + & + & + & + & - & - & + & - & + & + \\
+ & + & - & + & + & - & - & + & + & + & - & + & + & - & - & + & + & - & + & - & + & - \\
+ & - & - & - & - & + & - & - & - & + & - & - & + & + & + & - & + & + & + & + & - & + \\
+ & - & + & + & - & + & + & - & + & + & - & - & - & + & + & - & + & - & - & + & + & + \\
- & + & + & - & + & + & - & - & - & - & + & - & + & + & + & - & + & - & + & + & + & + \\
+ & + & + & - & + & + & + & - & + & + & + & - & + & - & - & - & - & + & - & + & - & + \\
+ & - & + & + & - & - & + & + & - & - & + & + & + & + & + & + & - & - & + & - & + & - \\
- & - & - & + & + & + & + & + & + & + & - & - & - & + & + & + & + & - & + & - & + & - \\
- & + & + & + & - & - & + & + & + & - & + & + & + & - & + & - & - & + & + & - & + & - \\
- & - & + & - & - & - & - & - & - & - & + & - & + & - & + & - & - & - & + & - & - & - \\
- & + & + & - & + & - & - & - & + & + & - & + & - & - & + & - & - & + & + & + & - & + \\
\end{array}\right)\)}
\end{minipage}
\hfill
\begin{minipage}[t]{0.485\textwidth}
\centering
{\footnotesize $H_{22}=$}\par\vspace{-0.2mm}
{\fontsize{5.8}{5.9}\selectfont
\setlength{\arraycolsep}{0.55pt}\renewcommand{\arraystretch}{0.68}
\(\left(\begin{array}{*{22}{c}}
+ & - & + & + & - & + & - & - & - & + & + & - & + & - & - & + & - & - & + & + & + & - \\
- & - & - & + & + & + & + & + & - & + & + & - & + & - & + & + & + & + & - & + & - & - \\
- & - & + & + & - & - & - & + & + & - & + & + & + & - & - & + & + & + & + & - & - & - \\
+ & + & + & + & - & + & - & + & - & + & + & + & - & - & + & - & - & - & + & - & + & + \\
- & - & - & - & + & + & + & + & - & - & + & + & - & + & + & - & + & - & - & + & + & - \\
+ & + & - & + & - & - & + & + & - & + & + & + & - & - & + & + & + & + & - & - & + & - \\
+ & - & + & + & - & - & - & - & + & + & - & + & - & + & - & + & - & + & + & + & + & + \\
- & + & - & - & + & + & - & - & + & + & + & - & + & + & - & + & + & + & + & + & - & - \\
- & + & - & + & + & - & + & - & + & + & + & + & + & + & - & - & + & - & - & - & - & - \\
- & + & + & - & - & + & - & + & + & - & + & + & + & + & - & - & - & + & + & + & + & - \\
- & + & - & + & - & + & + & + & - & + & - & + & + & + & - & + & + & + & + & + & + & + \\
+ & + & + & - & - & + & - & + & + & - & + & + & + & + & - & - & - & + & + & - & - & + \\
+ & + & - & + & + & + & + & - & - & - & - & + & - & + & + & + & + & - & - & + & - & + \\
- & + & + & - & + & + & - & - & + & + & + & - & + & + & - & - & + & + & - & + & - & + \\
+ & + & - & - & - & - & + & - & - & - & + & - & - & + & + & + & - & + & + & + & + & - \\
+ & + & - & + & + & - & + & + & - & + & + & - & - & - & + & + & - & + & - & - & + & + \\
+ & - & + & + & - & + & + & - & - & - & - & + & - & + & + & + & - & + & - & + & + & + \\
+ & + & + & + & - & + & + & + & - & + & + & + & - & + & - & - & - & - & + & - & + & - \\
- & + & - & + & + & - & - & + & + & - & - & + & + & + & + & + & + & - & - & + & - & + \\
- & - & - & - & + & + & + & + & + & + & + & - & - & - & + & + & + & + & - & + & - & + \\
- & - & + & + & + & - & - & + & + & + & - & + & + & + & - & + & - & - & + & + & - & + \\
- & - & - & + & - & - & - & - & - & - & - & + & - & + & - & + & - & - & - & + & - & - \\
\end{array}\right)\)}
\end{minipage}
\restoregeometry

The array determined by these slices satisfies the defining orthogonality conditions in all three coordinate directions, and hence is a Hadamard cube of order $22$.

\paragraph{Remark.} Since an ordinary Hadamard matrix of order $22$ does not exist, the example cannot be a $(22,3,1)$-Hadamard hypercube; in this sense the construction is genuinely of class $2$.

\section*{Disclosure on the use of generative AI}
The computational search that led to the example displayed in this note was carried out with assistance from ChatGPT using the GPT-5.6 Sol model (OpenAI). ChatGPT was also used to assist with formatting the explicit array and preparing an initial draft of this short note. The displayed array was checked directly against the defining orthogonality conditions.

\end{document}